\documentclass[12pt]{amsart}
\usepackage[utf8]{inputenc}
\usepackage{amsfonts,amsmath,amssymb,amscd,amsthm,url,enumerate,enumitem}
\usepackage{graphicx,epsf,afterpage}
\usepackage[english]{babel}

\usepackage{hyperref}
\hypersetup{
   colorlinks   = true, 
   urlcolor     = blue, 
   linkcolor    = blue, 
   citecolor   = red 
}

\theoremstyle{plain}
\newtheorem{theorem}{Theorem}[section]
\newtheorem{lemma}[theorem]{Lemma}
\newtheorem{proposition}[theorem]{Proposition}

\newtheorem{Kolmogorov}[theorem]{Kolmogorov's Theorem}
\newtheorem{Approx}[theorem]{Approximative Kolmogorov's Theorem}
\newtheorem{CountL}[theorem]{Countable Approximative Kolmogorov's Lemma}

\theoremstyle{definition}
\newtheorem{remark}[theorem]{Remark}

\newcommand\norm[1]{\ensuremath{\left\lVert#1\right\rVert}}
\newcommand\abs[1]{\ensuremath{\left\lvert#1\right\rvert}}
\newcommand\babs[1]{\ensuremath{\bigl\lvert#1\bigr\rvert}}

\def\R{{\mathbb R}} \def\Z{{\mathbb Z}}  
\def\tofive{1,\ldots,5}
\def\tofiveset{\{1,\ldots,5\}}

\begin{document}

\title{A structured proof of Kolmogorov's Superposition Theorem}
\author{S. Dzhenzher and A. Skopenkov}
\date{}

\thanks{Both authors: Moscow Institute of Physics and Technology. A. Skopenkov: Independent University of Moscow. 
\newline 
We would like to thank S. Shaposhnikov and anonymous reviewer for helpful discussions.
\newline Research supported by the Russian Foundation for Basic Research Grant No. 19-01-00169.}

\maketitle

\begin{abstract}
We present a well-structured detailed exposition of a well-known proof of the following celebrated result solving Hilbert's 13th problem on superpositions.
For functions of $2$ variables the statement is as follows.     
    
\textbf{Kolmogorov's Theorem.}
\emph{There are continuous functions $\varphi_1,\ldots,\varphi_5 \colon [\,0, 1\,]\to [\,0,1\,]$ such that for any continuous function $f\colon [\,0,1\,]^2\to\R$  there is a continuous function $h\colon [\,0,3\,]\to\R$ such that for any $x,y\in [\,0, 1\,]$ we have
\[
    f(x,y)=\sum\limits_{k=1}^5 h\left(\varphi_k(x)+\sqrt{2}\,\varphi_k(y)\right).
\]
}
The proof is accessible to non-specialists, in particular, to students familiar with only basic properties of continuous functions.
\end{abstract}

\section{Introduction}

The reader understand that a polynomial in some variables can be expressed using addition and multiplication (and constants).
Probably, the reader is familiar with expression of one Boolean function through another (see, e.~g., \cite{SZ}; however, this familiarity is not necessary to understand this article).

Kolmogorov's Theorem~\ref{thKolm} shows that any continuous function of two and more variables `can be expressed' using addition and continuous functions of one variable.
(A rigorous formulation of the concept `can be expressed' is not necessary for this theorem; however, the reader can find this formulation and its discussion, e.~g., in~\cite{Ar}.)
This theorem solves Hilbert's 13th problem.
We present a well-structured detailed exposition of the well-known proof of this important theorem.


\begin{Kolmogorov}\label{thKolm}
For any integer $n > 1$
    
there are real numbers $\alpha_1, \ldots, \alpha_n$ and continuous functions $\varphi_1,\varphi_2,\ldots,\varphi_{2n+1} \colon [\,0,1\,]\to [\,0,1\,]$ such that
    
for any continuous function $f\colon[\,0,1\,]^n\to\R$
    
there is a continuous function $h\colon\R\to\R$ such that
    \[
        f(x_1, \ldots, x_n) = \sum_{k=1}^{2n+1} h\left(\sum_{i=1}^n \alpha_i \varphi_k(x_i)\right).
    \]
\end{Kolmogorov}

\begin{remark}\label{r:hist}\leavevmode
(a) We can take $\alpha_1,\ldots,\alpha_n$ to be square roots of pairwise different prime numbers.
Or, for $n=2$, take $\alpha_1=1$ and $\alpha_2=\sqrt2$.  
E.~g., for the function $f(x, y) = x + \sqrt{2}\,y + 2$ we may take $h(x) = x$, $\varphi_1(x) = x$, $\varphi_2(x) = \frac{2}{1 + \sqrt{2}}$,  and $\varphi_3 = \varphi_4 =\varphi_5 =0$.
However, we do not know explicit $\varphi,h$ even for as simple functions as addition and multiplication.
And of course we do not know explicit universal $\varphi$.
    
(b) For classical expositions of a proof, for statement and discussion of Hilbert's 13th problem see e.~g. \cite{Ar, He, Ko}, \cite[Chapter 17]{LGM}, \cite[\S\S1--4]{St}.
In order to make our exposition well-structured we e.~g. explicitly introduce the notion of $\lambda$-prekolmogorov maps and explicitly state Approximative Kolmogorov's Theorem~\ref{thKolmApprox}. 
In order to make our exposition more detailed we e.~g. explicitly state  Propositions~\ref{propQuarters} and~\ref{propRationalSeparatingMap}, and present their simple proof.
    Cf.~\cite{BCM}.
    
(c) This topic is actively studied not only in analysis but also in topology and computer science, see e.~g. surveys \cite{St, Vi, Sk, Br, SH} and the references therein.
    
\end{remark}

\textbf{Outline of the proof of Kolmogorov's Theorem~\ref{thKolm}.}

We present a proof for $n=2$, taking $\alpha_1=1$ and $\alpha_2=\sqrt2$.
Proof for the general case is analogous. 
For a proof we need some conventions, notation and definitions.

\emph{In this text `function/map' means `continuous function/map'.}

Denote $I:=[\,0, 1\,]$.
We regard an ordered set of functions $\varphi_1,\varphi_2,\ldots,\varphi_5 \colon I\to I$ as a map (=vector-function) $\varphi\colon I\to I^5$. 
For a continuous function $h\colon [\,0,3\,]\to\R$ denote 
\[
    S_{\varphi} h (x,y) := \sum_{k=1}^5 h\left(\varphi_k(x)+\sqrt2\varphi_k(y)\right).
\] 

For a compact subset $M \subset \mathbb R^s$ and a function $f\colon M\to\mathbb R$ denote $\norm{f}:=\sup_{z\in M}\lvert f(z)\rvert$.\footnote{All functions further in the text have a compact domain, which is taken as $M$ for each of the functions (i.~e., for different functions the sets $M$ are different).} 

A \textbf{$\lambda$-prekolmogorov map} for a non-zero function $f\colon I^2\to\mathbb R$ is a map $\varphi\colon I\to I^5$ for which there exists a function $h\colon [\,0,3\,]\to\mathbb R$ such that
\[
    \norm{f - S_{\varphi}h} < \lambda \norm{f} \quad\text{and}\quad \norm{h} \leqslant \norm{f}.
\]
For the function $f\equiv0$ any map $\varphi\colon I\to I^5$ is $\lambda$-prekolmogorov.

\begin{Approx}\label{thKolmApprox}
There is a map $\varphi\colon I\to I^5$ which is $(7/8)$-prekolmogorov for any function $f\colon I^2\to\R$. 
\end{Approx}

\begin{proof}[Deduction of Theorem~\ref{thKolm} for $n=2$, $\alpha_1=1$ and $\alpha_2=\sqrt2$ from Theorem~\ref{thKolmApprox}]
    For non-ne\-gative integers $m$ define functions $h_m\colon[\,0,3\,]\to\R$ inductively.\footnote{Note that it is sufficient to define $h_m$ and $h$ on $[\,0,1+\sqrt2\,]$, but instead we write $[\,0,3\,]$ for brevity.}
    Let $h_0\equiv0$.  
    Let $h_m$ be a function obtained by applying Theorem \ref{thKolmApprox} to $f-\sum\limits_{k=0}^{m-1} S_\varphi h_k$.
    Then denoting $\lambda = 7/8$ we have 
    \begin{equation}
        \norm{f - \sum\limits_{k=0}^m S_\varphi h_k}
        < \lambda\norm{f - \sum\limits_{k=0}^{m-1} S_\varphi h_k} < \ldots < \lambda^m\norm{f}.
        \tag{*}\label{leqNorm}
    \end{equation}
    Since $\norm{h_m} \leqslant \norm{f - \sum\limits_{k=0}^{m-1} S_\varphi h_k} < \lambda^{m-1}\norm{f}$, the functional series $\sum\limits_{m=0}^\infty h_m$ uniformly converges to a function $h\colon[\,0,3\,]\to\R$.
    Then $\sum\limits_{m=0}^\infty S_\varphi h_m = S_\varphi\sum\limits_{m=0}^\infty h_m=S_\varphi h$.
    Since $\lambda^m\xrightarrow{m \to \infty} 0$, passing to the limit in \eqref{leqNorm} we obtain $\norm{f - S_{\varphi}h}= 0$, i.~e., $S_{\varphi}h=f$.
    Now extend $h$ to $\R$.
\end{proof}

\begin{CountL}\label{l:polyn}
    There is a map $\varphi\colon I\to I^5$, which is $(6/7)$-prekolmogorov for any polinomial $g\colon I^2\to\R$ with rational coefficients.
\end{CountL}

\begin{proof}[Proof of Theorem~\ref{thKolmApprox} modulo Lemma~\ref{l:polyn}]
    Take any functi\-on $f\colon I^2\to\R$ and the map $\varphi\colon I\to I^5$, given by Lemma~\ref{l:polyn}.  
    If $f \equiv0$, then $\varphi$ is $(7/8)$-prekolmogorov for $f$ by definition.
    Suppose now that $\norm{f} > 0$.
    Then for any function $g\colon I^2\to \R$ close enough to $\frac{111}{112}f$ we have $\norm{g} < \norm{f}$ and $\norm{f-g} < \frac{1}{56} \norm{f}$.
    The Weierstrass approximation theorem asserts that there is a polinomial $g\colon I^2\to\R$ with rational coefficients such that the inequalities above hold. 
    Since $\varphi$ is $(6/7)$-prekolmogorov for $g$, it follows that there is a function 
    \[
        h\colon[\,0,3\,]\to\R \quad\text{such that}\quad \norm h\leqslant\norm g\leqslant\norm f \quad\text{and}\quad \norm{g-S_\varphi h} < \frac67\norm g \leqslant \frac67\norm f.
    \]
    Then
    \[
        \norm{f - S_\varphi h} \leqslant \norm{f - g} + \norm{g - S_\varphi h} < \frac{1}{56}\norm{f} + \frac{6}{7}\norm{f} = \frac{7}{8}\norm{f}.
    \]
    Hence, $\varphi$ is $(7/8)$-prekolmogorov for $f$.
\end{proof}

\begin{lemma}\label{lemStability}
    For any function $f\colon I^2\to\R$ 

(stability) any map close enough to a $(6/7)$-prekolmogorov map for $f$, is a  $(6/7)$-prekol\-mogorov map for $f$; 

(approximation) for any map $\varphi\colon I\to I^5$ there is an arbitrary close $(6/7)$-prekolmogorov map for $f$.
\end{lemma}

For a function $f\colon I^2\to\R$ denote by $PK(f)$ the family of all $(6/7)$-prekolmogorov maps for $f$.

\begin{proof}[Deduction of Lemma~\ref{l:polyn} from Lemma~\ref{lemStability}]\footnote{The deduction uses the Baire category theorem. If you are not familiar with this theorem, just skip the deduction or replace the application of this theorem to the application of Cantor's intersection theorem.} 
Denote by $Q$ the set of all polynomials $I^2\to\R$ with rational coefficients.
    By Lemma~\ref{lemStability} for any $g \in Q$ the set $PK(g)$ is open and dense in the space of all continuous maps $I\to I^5$.
    It is well-known that this space is complete.
    Then by the Baire category theorem \cite{KF}, \cite[\S2]{Sk'} the countable intersection $\bigcap_{g\in Q} PK(g)$ is not empty.
Any map $\varphi$ from this intersection is as required.
\end{proof}

\section{Proof of Lemma~\ref{lemStability}}\label{s:prostability}

\emph{\textbf{Proof of the stability in Lemma~\ref{lemStability}.}} 
    If $f \equiv 0$, the statement is trivial. 
    So assume that $\norm{f}>0$. 
        
    Suppose $\psi\in PK(f)$ is a map.
    Then there is a function $h\colon[\,0,3\,]\to\R$ such that $\norm{f - S_\psi h} < \frac{6}{7}\norm{f}$. 
    Denote
    \[
        \varepsilon:= \frac{1}{5}\left(\frac{6}{7}\norm{f} - \norm{f - S_\psi h}\right) > 0.
    \]
    Since $h$ is uniformly continuous, there is $\delta$ such that $|h(x) - h(y)| < \varepsilon$ for $|x-y| < \delta$.
    It suffices to show that any map $\varphi$ which is $(\delta/3)$-close to $\psi$ is contained in $PK(f)$.
    
    For any $x,y\in I$ and $k = 1,\ldots,5$
    \begin{align*}
        &\left|\left(\psi_k(x) + \sqrt{2}\,\psi_k(y)\right) -
        \left(\varphi_k(x) + \sqrt{2}\,\varphi_k(y)\right)\right| <
        \left(1 + \sqrt{2}\right)\delta/3 < \delta.\\
    \intertext{Consequently,}
        &\left|h\left(\psi_k(x) + \sqrt{2}\,\psi_k(y)\right) -
        h\left(\varphi_k(x) + \sqrt{2}\,\varphi_k(y)\right) \right| < \varepsilon.
    \end{align*}
    Therefore $\norm{S_\psi h - S_\varphi h} < 5 \varepsilon$. Finally,
    \[
        \norm{f - S_\varphi h} \leqslant
        \norm{f - S_\psi h} + \norm{S_\psi h - S_\varphi h} <
        \frac{6}{7}\norm{f} - 5\varepsilon + 5\varepsilon =
        \frac{6}{7}\norm{f}.
    \]
    Hence $\varphi \in PK(f)$.
\hfill$\square$

\smallskip

\textbf{Preparation for the proof of the approximation in Lemma~\ref{lemStability}: rationally separating functions.}

Denote $[\,a,b\,] + c := [\,a+c,b+c\,]$ and $[\,a,b\,] \cdot d := [\,ad,bd\,]$.

\begin{figure}[ht]
    \centering
    \includegraphics[height=100mm]{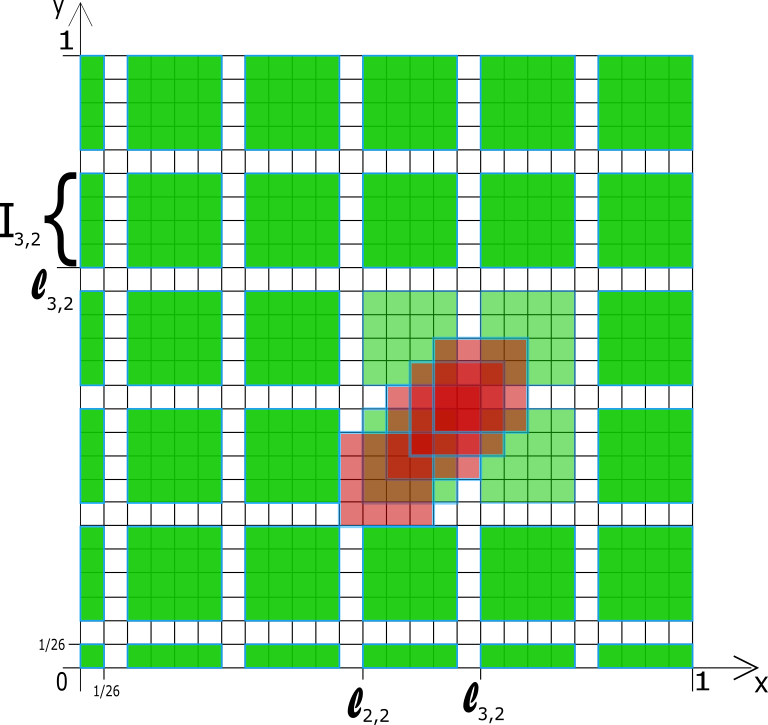}
    \caption{The green family of squares $I_{i,2}(26) \times I_{j,2}(26)$ parametrized by $i,j \in Z_2(26)$. The red squares are $I_{2,k} \times I_{1,k}$ for $k=1,3,4,5$.}
    \label{fig:squares}
\end{figure}

Take any $k = \tofive$. Denote
\[
    Z_k = Z_k(N) := \left[\,-1,\dfrac{N-k}5\,\right] \bigcap \Z.
\]
For every $j \in Z_k(N)$  denote
\[    
    I_{j,k} = I_{j,k}(N) := \dfrac{4I+5j + k}{N}, \quad l_{j,k} = l_{j,k}(N) := \dfrac{5j + k}{N} \quad\text{and}\quad l_{j,k}^0 := \max\{l_{j,k},0\}.
\]
Obviously, $l_{j,k}$ is the left end of the interval $I_{j,k}$.\footnote{For fixed $k$ and $N$ the segments $I_{j,k}(N)$ for any $j\in Z_k(N)$ have length $\frac4N$, they do not intersect and distances between neighboring segments are equal to $\frac1N$, and the first segments covers the point $0$, and the last one covers the point $1$.}

\begin{proposition}\label{propQuarters}
    Let $N$ be a positive integer. For any $(x, y) \in I^2$ and for at least three $k=\tofive$ there are $i, j \in Z_k$ such that
    \[
        (x, y) \in I_{i,k}(N) \times I_{j,k}(N).
    \]
\end{proposition}

\begin{proof}
    We have 
    \[
        (x, y) \in I_{i,k}(N) \times I_{j,k}(N) \quad\Leftrightarrow\quad 
        (Nx, Ny) \in I_{i,k}(1) \times I_{j,k}(1).
    \]
    Denote $m := \left[\frac{Nx}{5}\right]$ and denote by $r$ the remainder of the division of $[Nx]$ by $5$.
    Then $Nx \in [\,5m + r, 5m + r + 1)$ and for any $s = 0,1,2,3$
    we have
    \[
        Nx \in [\,5m + (r - s), 5m + (r - s + 4)\,] = 4I + 5m + (r-s).
    \]
    Therefore for $k \equiv r - s \pmod 5$ there is $i \in Z_k(N)$ such that $Nx \in I_{i,k}(1)$.
    Then for at least four $k = \tofive$ there is the appropriate $i$.
    Similarily to this, for at least four $k=\tofive$ there is $j\in Z_k(N)$ such that $Ny \in I_{j,k}(1)$.
    Then for at least three $k=\tofive$ there are $i,j \in Z_k(N)$ such that $(Nx, Ny) \in I_{i,k}(1) \times I_{j,k}(1)$.
\end{proof}

A function $\varphi\colon I\to I$ \emph{rationally separates} a number of pairwise disjoint closed intervals on the line if $\varphi$ has constant rational pairwise different values on the intervals.\footnote{Recall that in this text `function' means `continuous function', and `map' means `continuous map'.}

For a map $\varphi\colon I\to I^5$ denote $\norm{\varphi} := \max\limits_{k=\tofive} \norm{\varphi_k}$.

\begin{proposition}\label{propRationalSeparatingMap}
    Let $\varepsilon$ be a positive number.
    
    \begin{enumerate}[ref=\text{\alph*}, label=\text{(\alph*)}]
        \item\label{pointRationalSeparating}
        Let $\psi\colon I\to I$ be a function and $k\in\tofiveset$.
        
        Then there is an integer $N_0 = N_0(\psi, k) > 0$ such that for any integer $N>N_0$ and for any finite set $G$ of rational numbers there is a function $\varphi\colon I\to I$ such that
        
        \begin{itemize}
            \item $\norm{\varphi - \psi} < \varepsilon$;
            
            \item $\varphi$ rationally separates the family of intervals $I_{j,k} = I_{j,k}(N)$ parametrized by $j\in Z_k(N)$;
            
            \item $\varphi(I_{j,k}) \not\in G$ for any $j \in Z_k(N)$.
        \end{itemize}
        
        \item\label{pointRationalSeparatingMap}
        Let $\psi\colon I\to I^5$ be a map. 
        
        Then there is an integer $N_0 > 0$ such that for any integer $N>N_0$ there is a map $\varphi\colon I\to I^5$ such that
        \begin{itemize}
            \item $\norm{\varphi - \psi} < \varepsilon$;
            
            \item for every $k=\tofive$ the function $\varphi_k$ rationally separates the family of intervals $I_{j,k} = I_{j,k}(N)$ parametrized by $j\in Z_k(N)$;
            
            \item numbers $\varphi_k(I_{j,k})$ are pairwise different for different $k\in\tofiveset$ and $j \in Z_k(N)$.
        \end{itemize}
    \end{enumerate}
        
\end{proposition}

\begin{proof}[Proof of part~\eqref{pointRationalSeparating}.]
    Since $\psi$ is uniformly continuous, there is $\delta > 0$ such that $\abs{\psi(x) - \psi(y)} < \dfrac{\varepsilon}{2}$ for $\abs{x - y} < \delta$.
    Define $N_0: = \left\lceil\dfrac{5}{\delta}\right\rceil$.
    Take any $N>N_0$.
    Take a piecewise linear function $\varphi\colon I\to I$ such that
    \begin{itemize}
        \item for each interval $I_{j,k}$ the value $\varphi(I_{j,k})$ is rational and is $\frac{\varepsilon}{2}$-close to $\psi\!\left(l_{j,k}^0\right)$;
        
        \item all values $\varphi(I_{j,k})$ are pairwise different and do not lie in $G$;
        
        \item $\varphi$ is linear on the gaps between the intervals.
    \end{itemize}
    
    Since $0 \in I_{-1,k}$ and $1 \in I_{\max Z_k, k}$ for any $k=\tofive$, the definition of $\varphi(0)$ and $\varphi(1)$ is meaningful.
    
    Now $\varphi$ rationally separates the family of intervals $I_{j,k}$ and its values on the intervals are not from $G$. It remains to show that
    \begin{equation}
        \abs{\varphi(x) - \psi(x)} < \varepsilon\quad \text{for any}\quad x \in I. \tag{**}\label{ineqEpsilonClose}
    \end{equation}
    
    For any point $x \in I_{j,k}$, since $\varphi(x) = \varphi(l_{j,k}^0)$ and $\abs{x - l_{j,k}^0} \leqslant \frac{4}{N} < \delta$, we have
    \[
        \bigl|\varphi(x) - \psi(x)\bigr| \leqslant
        \bigl|\varphi(l_{j,k}^0) - \psi(l_{j,k}^0)\bigr| +
        \bigr|\psi(l_{j,k}^0) - \psi(x)\bigr| <
        \frac{\varepsilon}{2} + \frac{\varepsilon}{2} = \varepsilon.
    \]
    For any point $x$ which does not belong to any interval $I_{j,k}$ choose $j \in Z_k(N)$ such that $x$ lies between two intervals $I_{j,k}$ and $I_{j+1,k}$, i.~e., $ l_{j,k}+\frac4N < x < l_{j+1,k} = l_{j+1,k}^0$.  
    Ends of intervals are related by the obvious equation $l_{j,k} + \frac{5}{N} = l_{j+1,k}$.
    Then there is $\alpha \in (0,1)$ such that $x=l_{j+1,k}-\frac{\alpha}N$. 
    Now~\eqref{ineqEpsilonClose} follows because
    \[
    \begin{split}
        \abs{\varphi(x) - \psi(x)} &\leqslant
        \alpha\abs{\varphi\left(l_{j+1,k}-\tfrac1N\right) - \psi(x)} +
        (1 - \alpha)\abs{\varphi(l_{j+1,k}) - \psi(x)}
        < \\ &<
        \alpha \varepsilon + (1 - \alpha) \varepsilon
        = \varepsilon.
    \end{split}
    \]
    Here the first inequality follows because $\varphi$ is linear on the gap $l_{j+1,k}+[\,-\frac1N,0\,]$ (i.~e., $\varphi(x) = \alpha\varphi(l_{j+1,k}-\frac1N)+(1-\alpha)\varphi(l_{j+1,k})$).
    The second inequality if proved as follows.
    Since $\varphi(l_{j+1,k}-\frac1N) = \varphi(l_{j,k}^0)$ and $\abs{x - l_{j,k}^0} \leqslant \frac{5}{N} < \delta$, we have
    \[
        \abs{\varphi(l_{j+1,k}-\tfrac1N) - \psi(x)} \leqslant
        \abs{\varphi(l_{j,k}^0) - \psi(l_{j,k}^0)} + \abs{\psi(l_{j,k}^0) - \psi(x)} <
        \frac{\varepsilon}{2} + \frac{\varepsilon}{2} = \varepsilon.
    \]
    Since $\abs{x - l_{j+1,k}} \leqslant \frac{1}{N} < \delta$, we have
    \[
        \abs{\varphi(l_{j+1,k}) - \psi(x)} \leqslant
        \abs{\varphi(l_{j+1,k}) - \psi(l_{j+1,k})} + \abs{\psi(l_{j+1,k}) - \psi(x)} <
        \frac{\varepsilon}{2} + \frac{\varepsilon}{2} = \varepsilon.
    \]
\end{proof}

\begin{proof}[Proof of part~\eqref{pointRationalSeparatingMap}.]
    Choose $N_0 = \max\limits_{k=\tofive} N_0(\psi_k, k)$, where $N_0(\psi_k, k)$ come from applications of part~\eqref{pointRationalSeparating}.
    Take any $N > N_0$.
    Apply part~\eqref{pointRationalSeparating} to $\psi_1$, $k=1$ and $G = \varnothing$.
    We obtain a function $\varphi_1$. 
    Then apply part~\eqref{pointRationalSeparating} to $\psi_2$, $k=2$ and $G = \bigl\{\varphi_1(I_{j,1})\bigr\}_{j \in Z_1}$. 
    We obtain a function $\varphi_2$. 
    Analogously, on the $k$th step we apply part~\eqref{pointRationalSeparating} to the function $\psi_k$ and the set $G = \bigcup\limits_{s=1}^{k-1} \bigl\{\varphi_s(I_{j,s})\bigr\}_{j \in Z_s}$. 
    We obtain a function $\varphi_k$. 
    After five steps we obtain the required map $\varphi:=(\varphi_1,\ldots,\varphi_5)$.
\end{proof}

\emph{\textbf{Proof of the approximation in Lemma \ref{lemStability}.}}
    If $f \equiv 0$, the statement is trivial. 
    So assume that $\norm{f}>0$. 
    
    Fix any map $\psi\colon I\to I^5$ and any $\varepsilon > 0$.
    We should prove that there is $\varphi \in PK(f)$ which is $\varepsilon$-close to $\psi$.
    
    Apply Proposition~\ref{propRationalSeparatingMap}.\ref{pointRationalSeparatingMap} to the map $\psi$ and $\varepsilon$. 
    We obtain a number $N_0$.
    Since $f$ is uniformly continuous, there is an integer $N > N_0$ such that $\abs{f(x, y) - f(x', y')} < \frac{1}{6}\norm{f}$ for $\abs{x - x'} < \frac{4}{N}$ and $\abs{y - y'} < \frac{4}{N}$.
    Apply Proposition~\ref{propRationalSeparatingMap}.\ref{pointRationalSeparatingMap} to this $N$.
    We obtain a map $\varphi\colon I\to I^5$.
    
    For each $k$ define the function
    \[
        \tilde \varphi_k\colon I^2\to [\,0,3\,]\quad \text{by}\quad \tilde\varphi_k(x, y) = \varphi_k(x)+\sqrt{2}\,\varphi_k(y).
    \]

    In the next paragraph we prove that $\tilde\varphi_k$ has different constant values on different squares that are cartesian products of the intervals $I_{j,k}$.
    
    Indeed, suppose for some $x_1, x_2, y_1, y_2$ belonging to some of the intervals we have
    $\varphi_k(x_1) + \sqrt{2}\,\varphi_k(y_1) = \varphi_k(x_2) + \sqrt{2}\,\varphi_k(y_2)$. By definition of rational separability, the numbers $\varphi_k(x_1), \varphi_k(x_2), \varphi_k(y_1)$ and $\varphi_k(y_2)$ are rational.
    Then $\varphi_k(x_1) = \varphi_k(x_2)$ and $\varphi_k(y_1) = \varphi_k(y_2)$.
    Consequently, the pairs $x_1, x_2$ and $y_1, y_2$ are from the same intervals. Then the points $(x_1, y_1)$ and $(x_2, y_2)$ are from the same square.
    
    In the next paragraph we prove that the numbers $\tilde\varphi_k(I_{i,k} \times I_{j,k})$ are  pairwise different for different triples $(i,j,k)$. 

    Indeed, suppose $\tilde\varphi_k(I_{i,k} \times I_{j,k}) = \tilde\varphi_n(I_{p,n} \times I_{q,n})$. Then $\varphi_k(I_{i,k}) = \varphi_n(I_{p,n})$. So $k=n$ and $i=p$. Analogously $j=q$.
    
    Take a piecewise linear function $h\colon[\,0,3\,]\to\R$ such that $h\bigl(\tilde\varphi_k(I_{i,k} \times I_{j,k})\bigr) = \frac{1}{3}f(l_{i,k}^0, l_{j,k}^0)$ for any $k=\tofive$ and $i,j \in Z_k$.
    Then $\norm{h} \leqslant \frac{1}{3}\norm{f}$.
    Now the lemma follows because for any $z \in I^2$
    \[
        \babs{S_\varphi h(z) - f(z)} \leqslant \frac{1}{6}\norm{f} + \frac{2}{3}\norm{f} < \frac{6}{7}\norm{f}.
    \]
    Here the second inequality is obvious.
    The first inequality is proved as follows.
    Since $f$ is uniformly continuous, $\abs{f(z) - f(z')} < \frac{1}{6}\norm{f}$ for any $k=1,\ldots,5$, any $i,j \in Z_k$ and any $z,z' \in I_{i,k} \times I_{j,k}$.
    By Proposition~\ref{propQuarters}, there are at least three triples $(i_s,j_s,k_s) \in Z_{k_s} \times Z_{k_s} \times \tofiveset$, $s = 1,2,3$, such that $I_{i_s,k_s}\times I_{j_s,k_s} \ni z$.
    For them $h\bigl(\tilde\varphi_{k_s}(z)\bigr) = \dfrac{1}{3}f(l_{i_s,k_s}^0, l_{j_s,k_s}^0)$. These values differ from $\dfrac{1}{3}f(z)$ by at most  $\dfrac{1}{18}\norm{f}$.
    Hence
    \[
     \abs{\sum_{s=1}^3 h\bigl(\tilde\varphi_{k_s}(z)\bigr) - f\bigl(z\bigr)} \leqslant \dfrac{1}{6}\norm{f}.
    \]
    Another two values $h\bigl(\tilde\varphi_k(z)\bigr)$ do not exceed $\dfrac{1}{3}\norm{f}$ by the absolute value.
    \hfill$\square$


\begin{thebibliography}{99}



\bibitem[Ar]{Ar} \emph{V.I. Arnold},
On representation of continuous functions of several variables by superpositions of functions of less variables (in Russian). Mat. Prosveschenie, ser. 2, 1958.
\newline
\url{http://ilib.mccme.ru/djvu/mp2/mp2-3.htm}



\bibitem[BCM]{BCM}
\emph{A. Belov, A. Chilikov, I. Mitrofanov, S. Shaposhnikov and A. Skopenkov,}
13th Hilbert Problem on superpositions of functions,
\url{http://www.turgor.ru/lktg/2016/5/index.htm}

\bibitem[Br]{Br}
\emph{V. Brattka},
From Hilbert’s 13th Problem to the theory of neural networks: constructive aspects of Kolmogorov’s Superposition Theorem, Kolmogorov’s Heritage in Mathematics, pp 253-280. Berlin, 2007. 





\bibitem[He]{He} \emph{T. Hedberg}, The Kolmogorov superposition theorem, Appendix II to H.S.Shapiro, Topics in Approximation Theory. Lecture Notes in Math., 1971, V. 187, P. 267--275.


\bibitem[KF]{KF}
\emph{A.N. Kolmogorov and S.V. Fomin}, 
Elements of the Theory of Functions and Functional Analysis, Dover Publ., 1999.  
 
\bibitem[Ko]{Ko}
\emph{A.N. Kolmogorov},
On representation of continuous functions of several variables by superpositions of functions of less variables and addition (in Russian). Dokl. AN SSSR, 114:5 (1957), 953-965.

\bibitem[LGM]{LGM}
\emph{G.G. Lorentz, M.v. Golitschek, Y. Makovoz},
Constructive Approximation: Advanced Problems.  Springer,  New York, 1996. 
 


\bibitem[SH]{SH}
\emph{J. Schmidt-Hieber},
The Kolmogorov–Arnold representation theorem revisited. Neural Networks, 137 (2021), 119--126.

\bibitem[Sk]{Sk} \emph{A. Skopenkov},
Basic embeddings and Hilbert's 13th problem (in Russian), Mat. Prosveshenie, 14 (2010), 143--174. arXiv:1001.4011. 
Abridged English translation:
arXiv:1003.1586. 

\bibitem[Sk']{Sk'}
\emph{A. Skopenkov},
Ambient Homogeneity, MCCME, Moscow, 2012. arXiv:1003.5278. 
  
\bibitem[St]{St}
\emph{Y. Sternfeld}, Hilbert's 13th problem and dimension, Lect. Notes Math. 1376 (1989) 1--49.  

\bibitem[SZ]{SZ} \emph{M. Skopenkov and A. Zaslavsky (editors)}, 
Mathematics via Problems: Part 3: Combinatorics, AMS, MSRI Mathematical Circles Library, to appear. 

\bibitem[Vi]{Vi}
\emph{A.G. Vitushkin}, Hilbert's 13th problem and related questions, Russian Math. Surveys,~59:1 (2004),~11--25.

 
\end{thebibliography}
\end{document}


\title[\mytitle]{\mytitle}
\author{С. В. Дженжер и А. Б. Скопенков}
\date{}

\maketitle


\begin{abstract}
    Мы представляем хорошо структурированное детальное изложение известного доказательства важного результата, являющегося решением 13-й проблемы Гильберта о суперпозициях.
    Для функций двух переменных он формулируется так. 
        
    \textbf{Теорема Колмогорова.}
    \emph{Существуют непрерывные функции
    $\varphi_1,\ldots,\varphi_5 \colon [\,0, 1\,]\to [\,0,1\,]$
    такие, что для любой непрерывной функции
    $f\colon [\,0,1\,]^2\to\R$ 
    существует непрерывная функция
    $h\colon [\,0,3\,]\to\R$
    такая, что для любых $x,y\in [\,0, 1\,]$
    \[
        f(x,y)=\sum\limits_{k=1}^5 h\left(\varphi_k(x)+\sqrt{2}\,\varphi_k(y)\right).
    \]
    }
    Доказательство доступно неспециалистам, в частности, студентам, знакомым только с основными свойствами непрерывных функций.
\end{abstract}

\section{Введение}

Читатель понимает, что многочлен от нескольких переменных может быть выражен через сложение и умножение (с использованием констант). Возможно, читатель знаком с выражением одних функций алгебры логики через другие (см., например, \cite[\S24.5]{ZSS}; впрочем, этого знакомства не нужно для понимания настоящей статьи). 

Теорема Колмогорова \ref{thKolm} показывает, что любая непрерывная функция двух и более переменных <<может быть выражена>> через непрерывные функции одной переменной и сложение. 
(Четкая формулировка понятия <<может быть выражена>> не нужна для этой теоремы; впрочем, читатель может найти эту формулировку и ее обсуждение, например, в \cite{Ar}.) 
Эта теорема решает 13-ю проблему Гильберта. 
Мы представляем хорошо структурированное детальное изложение известного доказательства этой важной теоремы.\footnote{Оба автора: Московский Физико-Технический Институт. А. Б. Скопенков: Независимый Московский Университет. 
\newline 
Благодарим С.В. Шапошникова и анонимного рецензента за полезные обсуждения.
\newline
Поддержано грантом Российского Фонда Фундаментальных Исследований No. 19-01-00169.}
 
\begin{theorem}[Колмогоров]\label{thKolm}
    Для любого целого $n>1$ 
    
    существуют действительные числа $\alpha_1,\ldots,\alpha_n$ и
    непрерывные функции
    \linebreak
    $\varphi_1,\varphi_2,\ldots,\varphi_{2n+1} \colon [\,0,1\,]\to [\,0,1\,]$
    такие, что 
    
    для любой непрерывной функции $f\colon[\,0,1\,]^n\to\R$ 
    
    существует непрерывная функция $h\colon\R\to\R$ 
    такая, что 
    \[
        f(x_1, \ldots, x_n) = \sum_{k=1}^{2n+1} h\left(\sum_{i=1}^n \alpha_i \varphi_k(x_i)\right).
    \]
\end{theorem}

\begin{remark}\label{r:hist}\leavevmode
    (a) В качестве $\alpha_1,\ldots,\alpha_n$ мы можем взять квадратные корни из попарно различных простых чисел.
    Или, для $n=2$, взять $\alpha_1=1$ и $\alpha_2=\sqrt2$.  
    Например, для функции $f(x,y)=x+\sqrt{2}y+2$ мы можем взять $h(x)=x$, $\varphi_1(x)=x$, $\varphi_2(x)=\frac{2}{1+\sqrt{2}}$ и $\varphi_3=\varphi_4=\varphi_5=0$.
    Тем не менее, мы не знаем функций $\varphi$, $h$ даже для таких простых функций $f$, как сложение и умножение. 
    И, конечно, мы не знаем явного выражения для <<универсальных функций>> $\varphi_1,\ldots,\varphi_{2n+1}$.
        
    (b) Классическое изложение доказательства, формулировка и обсуждения 13-й проблемы Гильберта приведены, например, в \cite{Ar, He, Ko}, \cite[Chapter 17]{LGM}, \cite[\S1-\S4]{St}.
    Чтобы сделать наше изложение хорошо структурированным, мы, например, явно ввели понятие $\lambda$-предколмогоровского отображения и явно сформулировали аппроксимативную теорему Колмогорова \ref{thKolmApprox}. 
    Чтобы сделать наше изложение более детальным, мы, например, явно сформулировали утверждения \ref{propQuarters} и \ref{propRationalSeparatingMap}, и привели их простые доказательства.
    Ср. с \cite{BCM}.     
        
    (c) Эта тема активно изучается не только в анализе, но также и в топологии и в компьютерной науке, см., например, работы \cite{St, Vi, Sk, Br, SH} и ссылки в них.
    
    В 27-м и 28-м выпусках Матпросвещения опубликованы задачи соответственно 27.4 (с. 234) и 27.4' (с. 248), связанные по тематике с 13-й проблемой Гильберта.  
\end{remark}

\textbf{План доказательства теоремы \ref{thKolm} Колмогорова.}

Мы представляем доказательство для $n=2$, полагая $\alpha_1=1$ и $\alpha_2=\sqrt2$.
Доказательство в общем случае аналогично.
Для доказательства нам необходимы некоторые соглашения, обозначения и определения.

\emph{Далее в тексте <<функция/отображение>> означет <<непрерывная функция/отображение>>.}

Обозначим $I:=[\,0, 1\,]$.
Мы рассматриваем упорядоченный набор функций $\varphi_1,\varphi_2,\ldots,\varphi_5 \colon I\to I$ как отображение (=вектор-функцию) $\varphi\colon I\to I^5$.
Для функции $h\colon [\,0,3\,]\to\R$ обозначим 
\[
    S_{\varphi} h (x,y) := \sum_{k=1}^5 h\left(\varphi_k(x)+\sqrt2\varphi_k(y)\right).
\] 
Для компактного подмножества $M\subset\mathbb R^s$ и функции $f\colon M\to\mathbb R$ обозначим $\norm{f}:=\sup_{z\in M}|f(z)|$.\footnote{Все функции в дальнейшем имеют компактную область определения, которая и берётся в качестве $M$ для каждой из функций (т.~е. для разных функций множества $M$ разные).} 
Назовём \textbf{$\lambda$-предколмогоровским отображением} для ненулевой функции $f\colon I^2\to\R$ отображение $\varphi\colon I\to I^5$, для которого найдётся функция $h\colon [\,0,3\,]\to\R$ такая, что
\[
    \norm{f - S_{\varphi}h} < \lambda \norm{f} \quad\text{и}\quad \norm{h} \leqslant \norm{f}.
\]
Для функции $f\equiv0$ любое отображение $\varphi\colon I\to I^5$ является $\lambda$-предколмогоровским.

\begin{theorem}[Аппроксимативная теорема Колмогорова]\label{thKolmApprox}
Существует отображение $\varphi\colon I\to I^5$, являющееся $(7/8)$-предколмогоровским для любой функции $f\colon I^2\to\R$. 
\end{theorem}

\begin{proof}[Вывод теоремы \ref{thKolm} для $n=2$, $\alpha_1=1$ и $\alpha_2=\sqrt2$ из теоремы \ref{thKolmApprox}]
Для целых неотрицательных $m$ определим функции $h_m\colon[\,0,3\,]\to\R$ индуктивно.\footnote{Заметим, что нам достаточно определить $h_m$ и $h$ на $[\,0,1+\sqrt2\,]$, но мы вместо этого пишем $[\,0,3\,]$ для краткости.}
Положим $h_0\equiv0$.
Обозначим через $h_m$ функцию, полученную применением теоремы \ref{thKolmApprox} к $f-\sum\limits_{k=0}^{m-1} S_\varphi h_k$.
Тогда, обозначив $\lambda = 7/8$, имеем
\begin{equation}
    \norm{f - \sum\limits_{k=0}^m S_\varphi h_k}
    < \lambda\norm{f - \sum\limits_{k=0}^{m-1} S_\varphi h_k} < \ldots < \lambda^m\norm{f}.  \tag{*}\label{leqNorm}
\end{equation}
Так как $\norm{h_m} \leqslant \norm{f - \sum\limits_{k=0}^{m-1} S_\varphi h_k} < \lambda^{m-1}\norm{f}$, то функциональный ряд $\sum\limits_{m=0}^\infty h_m$ сходится равномерно к функции $h\colon[\,0,3\,]\to\mathbb R$.
Тогда $\sum\limits_{m=0}^\infty S_\varphi h_m = S_\varphi\sum\limits_{m=0}^\infty h_m=S_\varphi h$.
Так как $\lambda^m\xrightarrow{m \to \infty} 0$, то переходя к пределу в \eqref{leqNorm}, получаем $\norm{f - S_{\varphi}h}= 0$, то есть $S_{\varphi}h=f$. 
Продолжим функцию $h\colon[\,0,3\,]\to\mathbb R$ до функции $h\colon\mathbb R\to\mathbb R$.
\end{proof}

\begin{lemma}[Счетная аппроксимативная лемма  Колмогорова]\label{l:polyn}
    Существует отображение $\varphi\colon I\to I^5$, являющееся $(6/7)$-предколмогоровским для любого многочлена $g\colon I^2\to\R$ с рациональными коэффициентами. 
\end{lemma}

\begin{proof}[Доказательство теоремы \ref{thKolmApprox} по модулю леммы \ref{l:polyn}] 
    Возьмем любые функцию $f\colon I^2\to\mathbb R$ и отображение $\varphi\colon I\to I^5$, данное леммой \ref{l:polyn}.  
    Если $f\equiv0$, то $\varphi$ является $(7/8)$-предколмогоровским для $f$ по определению.
    Теперь предположим, что $\norm{f}>0$.
    Тогда для любых достаточно близких к $\frac{111}{112}f$ функций $g\colon I^2\to \R$  выполнено $\norm{g} < \norm{f}$ и $\norm{f-g} < \frac{1}{56} \norm{f}$. 
    Из теоремы Вейерштрасса об аппроксимации непрерывных функций многочленами следует, что существует многочлен $g\colon I^2\to\R$ с рациональными коэффициентами, удовлетворяющий этим неравенствам.
    Так как $\varphi$ является $(6/7)$-предколмогоровским для $g$, то найдётся функция 
    \[
        h\colon[\,0,3\,]\to\R \quad\text{такая, что}\quad \norm{h} \leqslant \norm{g} \leqslant \norm{f} \quad\text{и}\quad \norm{g - S_\varphi h} < \frac{6}{7} \norm{g} \leqslant \frac{6}{7}\norm{f}.
    \]
    Тогда
    \[
        \norm{f - S_\varphi h} \leqslant \norm{f - g} + \norm{g - S_\varphi h} < \frac{1}{56}\norm{f} + \frac{6}{7}\norm{f} = \frac{7}{8}\norm{f}.
    \]
    Следовательно, $\varphi$ является $(7/8)$-предколмогоровским для $f$.
\end{proof}

\begin{lemma}\label{lemStability}
    Для любой функции $f\colon I^2\to\R$ 
    
    (стабильность) любое отображение, достаточно близкое к $(6/7)$-предколмогоровскому для $f$, само является $(6/7)$-предколмогоровским для $f$; 
    
    (аппроксимация) для любого отображения  $\varphi\colon I\to I^5$ найдется сколь угодно близкое $(6/7)$-предколмогоровское отображение для $f$. 
\end{lemma}

Для функции $f\colon I^2\to\R$ обозначим через $PK(f)$ семейство всех $(6/7)$-предколмогоровских отображений для $f$.  

\begin{proof}[Вывод леммы \ref{l:polyn} из леммы \ref{lemStability}]\footnote{Вывод использует  теорему Бэра о категории. Если вы не знакомы с этой теоремой, то просто пропустите этот вывод или замените применение этой теоремы на применение принципа вложенных шаров.}
    Обозначим через $Q$ семейство всех многочленов $I^2\to\R$ с рациональными коэффициентами.
    По лемме \ref{lemStability}, для любой функции $g\in Q$ множество $PK(g)$ открыто и всюду плотно в пространстве всех непрерывных отображений $I\to I^5$.
    Известно, что это пространство полно.
    Тогда по теореме Бэра о категории \cite[\S2.3.3]{KF}, \cite[\S2]{Sk'}, счётное пересечение $\bigcap_{g\in Q} PK(g)$ непусто.
    Любое отображение $\varphi$ из этого пересечения является искомым.
\end{proof}

\section{Доказательство леммы  \ref{lemStability}}\label{s:prostability}

\emph{\textbf{Доказательство стабильности в лемме \ref{lemStability}.}}
    Если $f \equiv 0$, то утверждение тривиально.
    Поэтому предположим, что $\norm{f}>0$.
            
    Предположим, $\psi\in PK(f)$ "--- некоторое отображение.
    Тогда существует функция $h\colon[\,0,3\,]\to\R$ такая, что $\norm{f - S_\psi h} < \frac{6}{7}\norm{f}$. 
    Обозначим
    \[
        \varepsilon := \frac{1}{5}\left(\frac{6}{7}\norm{f} - \norm{f - S_\psi h}\right) > 0.
    \]
    Так как $h$ равномерно непрерывна, то найдётся такое $\delta$, что $|h(x) - h(y)| < \varepsilon$ при $|x-y| < \delta$.
    Достаточно показать, что любое отображение $\varphi$, являющееся $(\delta/3)$-близким к $\psi$, лежит в $PK(f)$.
        
    Для любых $x,y\in I$ и $k = 1,\ldots,5$
    \begin{align*}
        &\left|\left(\psi_k(x) + \sqrt{2}\,\psi_k(y)\right) -
        \left(\varphi_k(x) + \sqrt{2}\,\varphi_k(y)\right)\right| <
        \left(1 + \sqrt{2}\right)\delta/3 < \delta.\\
    \intertext{Следовательно,}
        &\left|h\left(\psi_k(x) + \sqrt{2}\,\psi_k(y)\right) -
        h\left(\varphi_k(x) + \sqrt{2}\,\varphi_k(y)\right) \right| < \varepsilon.
    \end{align*}
    Поэтому $\norm{S_\psi h - S_\varphi h} < 5 \varepsilon$. Наконец,
    \[
        \norm{f - S_\varphi h} \leqslant
        \norm{f - S_\psi h} + \norm{S_\psi h - S_\varphi h} <
        \frac{6}{7}\norm{f} - 5\varepsilon + 5\varepsilon =
        \frac{6}{7}\norm{f}.
    \]
    Следовательно, $\varphi \in PK(f)$.
\hfill$\square$

\smallskip

{\bf Подготовка к доказательству аппроксимации в лемме \ref{lemStability}: рационально разделяющие функции.}

Обозначим $[\,a,b\,] + c := [\,a+c,b+c\,]$ и $[\,a,b\,] \cdot d := [\,ad,bd\,]$.

\begin{figure}[ht]
\centering
\includegraphics[height=100mm]{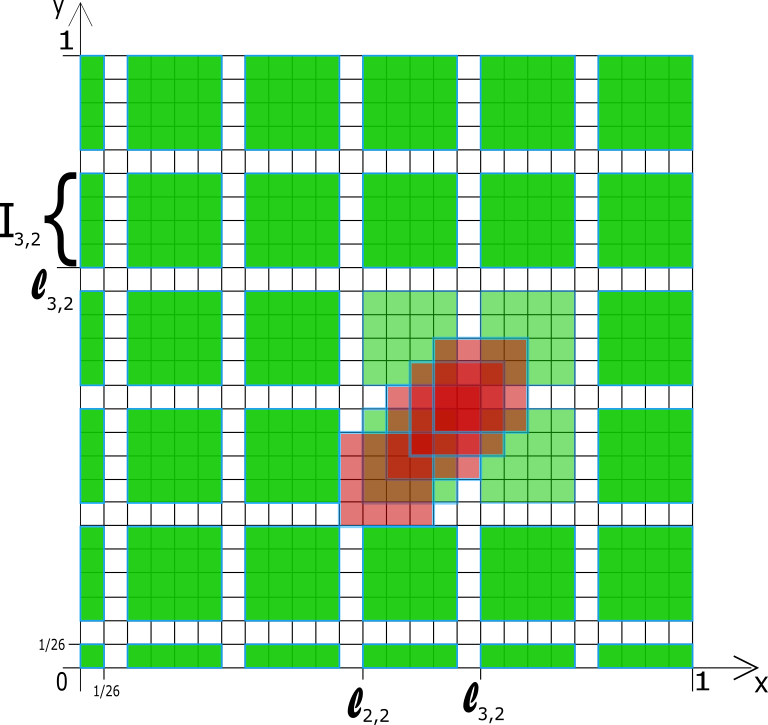}
\caption{Зелёные квадраты $I_{i,2}(26) \times I_{j,2}(26)$ параметризованы парами $i,j \in Z_2(26)$. Красные квадраты суть $I_{2,k} \times I_{1,k}$ для $k=1,3,4,5$.}
\label{fig:squares}
\end{figure}

Возьмём любое $k = \tofive$. Обозначим
\[
    Z_k = Z_k(N) := \left[\,-1,\dfrac{N-k}5\,\right] \bigcap \Z.
\]
Для каждого $j \in Z_k(N)$ обозначим
\[    
    I_{j,k} = I_{j,k}(N) := \dfrac{4I+5j + k}{N},\quad
    l_{j,k} = l_{j,k}(N) := \dfrac{5j + k}{N}
    \quad\text{и}\quad
    l_{j,k}^0 := \max\{l_{j,k},0\}.
\]
Очевидно, $l_{j,k}$ является левым концом отрезка $I_{j,k}$.\footnote{Для фиксированных $k$ и $N$ отрезки $I_{j,k}(N)$ при всевозможных $j\in Z_k(N)$ имеют длины $\frac4N$, они не пересекаются и расстояния между соседними отрезками равны $\frac1N$, причём самый первый отрезок покрывает точку $0$, а самый последний "--- точку $1$.}

\begin{proposition}\label{propQuarters}
    Пусть $N$ "--- целое положительное число. Для любых $(x, y) \in I^2$ и не менее чем  для трёх $k=\tofive$ найдутся $i, j \in Z_k(N)$ такие, что
    \[
        (x, y) \in I_{i,k}(N) \times I_{j,k}(N).
    \]
\end{proposition}

\begin{proof}
    Имеем
    \[
        (x, y) \in I_{i,k}(N) \times I_{j,k}(N) \quad\Leftrightarrow\quad 
        (Nx, Ny) \in I_{i,k}(1) \times I_{j,k}(1).
    \]
    Обозначим $m := \left[\frac{Nx}{5}\right]$ и обозначим через $r$ остаток от деления числа $[Nx]$ на $5$.
    Тогда $Nx \in [\,5m + r, 5m + r + 1)$, и для любого $s = 0,1,2,3$ имеем
    \[
        Nx \in [\,5m + (r - s), 5m + (r - s + 4)\,] = 4I + 5m + (r-s).
    \]
    Поэтому для $k \equiv r - s \pmod 5$ найдётся $i\in Z_k(N)$ такое, что $Nx \in I_{i,k}(1)$.
    Тогда для не менее чем чётырех $k = \tofive$ существует подходящее $i$. 
    Аналогично, для не менее чем четырёх $k=\tofive$ найдётся $j\in Z_k(N)$ такое, что $Ny \in I_{j,k}(1)$.
    Тогда для не менее чем трёх $k=\tofive$ найдутся такие $i,j \in Z_k(N)$, что $(Nx, Ny) \in I_{i,k}(1) \times I_{j,k}(1)$.
\end{proof}
  
Функция $\varphi\colon I\to I$ \emph{рационально разделяет} семейство попарно непересекающихся отрезков на прямой, если $\varphi$ принимает постоянные рациональные попарно различные значения на этих отрезках.\footnote{Напомним, что в данном тексте <<функция>> означает <<непрерывная функция>>, а <<отображение>> означает <<непрерывное отображение>>.} 
Для отображения $\varphi\colon I\to I^5$ обозначим $\norm{\varphi} = \max\limits_{k=\tofive} \norm{\varphi_k}$.

\begin{proposition}\label{propRationalSeparatingMap}
    Пусть $\varepsilon$ "--- некоторое положительное число.
    
    \begin{enumerate}[ref=\text{\alph*}, label=\text{(\alph*)}]
        \item\label{pointRationalSeparating}
        Пусть $\psi\colon I\to I$ "--- некоторая функция, и $k\in\tofiveset$.
            
        Тогда найдётся целое $N_0 = N_0(\psi, k) > 0$ такое, что для любого целого $N>N_0$ и для любого конечного множества $G$ рациональных чисел существует функция $\varphi\colon I\to I$ такая, что
            
        \begin{itemize}
            \item $\norm{\varphi - \psi} < \varepsilon$;
                    
            \item $\varphi$ рационально разделяет семейство отрезков $I_{j,k} = I_{j,k}(N)$ с параметром $j\in Z_k(N)$;
                    
            \item $\varphi(I_{j,k}) \not\in G$  для любого $j\in Z_k(N)$.
        \end{itemize}
            
        \item\label{pointRationalSeparatingMap}
        Пусть $\psi\colon I\to I^5$ "--- некоторое отображение.
            
        Тогда найдётся такое целое $N_0 > 0$, что для любого целого $N>N_0$ существует отображение $\varphi\colon I\to I^5$ такое, что
        \begin{itemize}
            \item $\norm{\varphi - \psi} < \varepsilon$;
                    
            \item для каждого $k=\tofive$ функция $\varphi_k$ рационально разделяет семейство отрезков $I_{j,k} = I_{j,k}(N)$ с параметром $j\in Z_k(N)$;
                    
            \item числа $\varphi_k(I_{j,k})$ попарно различны для различных $k\in\tofiveset$ и $j \in Z_k(N)$.
        \end{itemize}
    \end{enumerate}

\end{proposition}

\begin{proof}[Доказательство пункта \eqref{pointRationalSeparating}.]
    В силу равномерной непрерывности функции $\psi$, найдётся $\delta > 0$ такое, что $\abs{\psi(x) - \psi(y)} < \dfrac{\varepsilon}{2}$ при $\abs{x - y} < \delta$.
    Положим $N_0: = \left\lceil\dfrac{5}{\delta}\right\rceil$. 
    Рассмотрим любое $N>N_0$.
    Возьмём кусочно-линейную функцию $\varphi\colon I\to I$ такую, что
    \begin{itemize}
        \item на каждом отрезке $I_{j,k}$ значение $\varphi(I_{j,k})$ рационально и $\frac{\varepsilon}{2}$-близко к $\psi\!\left(l_{j,k}^0\right)$;
                
        \item все значения $\varphi(I_{j,k})$ различны и не лежат в $G$;
                
        \item $\varphi$ линейна на промежутках между отрезками.
    \end{itemize}
        
    Поскольку $0 \in I_{-1,k}$ и $1 \in I_{\max Z_k, k}$ для любого $k=\tofive$, то определение значений $\varphi(0)$ и $\varphi(1)$ осмысленно.
        
    Теперь $\varphi$ рационально разделяет семейство отрезков $I_{j,k}$, и её значения на отрезках не лежат в $G$. 
    Остаётся показать, что
    \begin{equation}
        \abs{\varphi(x) - \psi(x)} < \varepsilon\quad \text{для любого}\quad x \in I.
    \tag{**}\label{ineqEpsilonClose}
    \end{equation}
    
    Для любой точки $x \in I_{j,k}$, так как $\varphi(x) = \varphi(l_{j,k}^0)$ и $|x - l_{j,k}^0| \leqslant \frac{4}{N} < \delta$, то 
    \[
        \abs{\varphi(x) - \psi(x)} \leqslant
        \abs{\varphi(l_{j,k}^0) - \psi(l_{j,k}^0)} +
        \abs{\psi(l_{j,k}^0) - \psi(x)} <
        \frac{\varepsilon}{2} + \frac{\varepsilon}{2} = \varepsilon.
    \]
    Для любой точки $x$, не принадлежащей ни одному отрезку $I_{j,k}$, выберем $j \in Z_k(N)$ такое, что $x$ находится между двумя отрезками $I_{j,k}$ и $I_{j+1,k}$, т.~е., $l_{j,k}+\frac4N < x < l_{j+1,k} = l_{j+1,k}^0$.  
    Концы отрезков связаны между собой очевидным соотношением $l_{j,k} + \frac{5}{N} = l_{j+1,k}$.
    Тогда найдётся $\alpha \in (0,1)$ такое, что $x=l_{j+1,k}-\frac{\alpha}N$. 
    Теперь \eqref{ineqEpsilonClose} следует из того, что
    \[
        \begin{split}
            \abs{\varphi(x) - \psi(x)} &\leqslant
            \alpha\abs{\varphi\left(l_{j+1,k}-\tfrac1N\right) - \psi(x)} + 
            (1 - \alpha)\abs{\varphi(l_{j+1,k}) - \psi(x)}
                    < \\ &<
            \alpha\varepsilon + (1-\alpha)\varepsilon
            = \varepsilon.
        \end{split}
    \]
    Здесь первое неравенство следует из того, что $\varphi$ линейна на промежутке $l_{j+1,k}+[\,-\frac1N,0\,]$ (т.е. из $\varphi(x) = \alpha\varphi(l_{j+1,k}-\frac1N)+(1-\alpha)\varphi(l_{j+1,k})$).
    Второе неравенство доказывается так.
    Так как $\varphi(l_{j+1,k}-\frac1N) = \varphi(l_{j,k}^0)$ и $\abs{x - l_{j,k}^0} \leqslant \frac{5}{N} < \delta$, то
    \[
        \abs{\varphi(l_{j+1,k}-\tfrac1N) - \psi(x)} \leqslant
        \abs{\varphi(l_{j,k}^0) - \psi(l_{j,k}^0)} + \abs{\psi(l_{j,k}^0) - \psi(x)} <
        \frac{\varepsilon}{2}+\frac{\varepsilon}{2} = \varepsilon.
    \]
    Так как $\abs{x - l_{j+1,k}} \leqslant \frac{1}{N} < \delta$, имеем
    \[
            \abs{\varphi(l_{j+1,k}) - \psi(x)} \leqslant
            \abs{\varphi(l_{j+1,k}) - \psi(l_{j+1,k})} + \abs{\psi(l_{j+1,k}) - \psi(x)} <
            \frac{\varepsilon}{2} + \frac{\varepsilon}{2} = \varepsilon.
    \]
\end{proof}

\begin{proof}[Доказательство пункта~\eqref{pointRationalSeparatingMap}.]
    Выберем $N_0 = \max\limits_{k=\tofive} N_0(\psi_k, k)$, где $N_0(\psi_k, k)$ приходят из применения пункта~\eqref{pointRationalSeparating} утверждения.
    Возьмём любое $N > N_0$.
    Применим пункт~\eqref{pointRationalSeparating} к $\psi_1$, $k=1$ и $G = \varnothing$. Получим некоторую функцию $\varphi_1$. 
    Теперь применим пункт~\eqref{pointRationalSeparating} к $\psi_2$, $k=2$ и $G = \bigl\{\varphi_1(I_{j,1})\bigr\}_{j \in Z_1}$. 
    Получим функцию $\varphi_2$. 
    Аналогично, на $k$-ом шаге применяем пункт~\eqref{pointRationalSeparating} к функции $\psi_k$ и множеству $G = \bigcup\limits_{s=1}^{k-1} \bigl\{\varphi_s(I_{j,s})\bigr\}_{j \in Z_s}$. 
    Так получим функцию $\varphi_k$. 
    После пяти шагов получим требуемое отображение $\varphi:=(\varphi_1,\ldots,\varphi_5)$.
\end{proof}

\emph{\textbf{Доказательство аппроксимации в Лемме \ref{lemStability}.}}
    Если $f \equiv 0$, то утверждение тривиально.
    Поэтому предположим, что $\norm{f}>0$. 
    
    Зафиксируем любое отображение $\psi\colon I\to I^5$ и любое $\varepsilon > 0$. 
    Нужно указать отображение $\varphi \in PK(f)$, $\varepsilon$-близкое к $\psi$.
        
    Применим утверждение~\ref{propRationalSeparatingMap}.\ref{pointRationalSeparatingMap} к отображению $\psi$ и числу $\varepsilon$.   
    Получим число $N_0$.
    В силу равномерной непрерывности функции $f$, найдётся такое целое число $N > N_0$, что $\abs{f(x, y) - f(x', y')} < \frac{1}{6}\norm{f}$ при $\abs{x - x'} < \frac{4}{N}$ и $\abs{y - y'} < \frac{4}{N}$.
    Применим утверждение~\ref{propRationalSeparatingMap}.\ref{pointRationalSeparatingMap} к этому числу. 
    Получим отображение $\varphi\colon I\to I^5$.
        
    Для каждого $k \in \tofiveset$ определим функцию
    \[
        \tilde \varphi_k\colon I^2\to [\,0,3\,]\quad \text{формулой}\quad \tilde\varphi_k(x, y) = \varphi_k(x)+\sqrt{2}\,\varphi_k(y).
    \]
    В следующем абзаце мы покажем, что $\tilde\varphi_k$ принимает различные постоянные значения на различных квадратах, образованных декартовым произведением отрезков $I_{j,k}$.
        
    В самом деле, предположим, что для некоторых $x_1, x_2, y_1, y_2$, принадлежащих некоторым из этих отрезков, $\varphi_k(x_1) + \sqrt{2}\,\varphi_k(y_1) = \varphi_k(x_2) + \sqrt{2}\,\varphi_k(y_2)$. 
    По определению рациональной разделяемости, числа $\varphi_k(x_1), \varphi_k(x_2), \varphi_k(y_1)$ и $\varphi_k(y_2)$ рациональны.
    Тогда $\varphi_k(x_1) = \varphi_k(x_2)$ и $\varphi_k(y_1) = \varphi_k(y_2)$.
    Следовательно, пары $x_1, x_2$ и $y_1, y_2$ принадлежат одним и тем же отрезкам. 
    Тогда точки $(x_1, y_1)$ и $(x_2, y_2)$ принадлежат одному квадрату.
        
    В следующем абзаце мы покажем, что числа $\tilde\varphi_k(I_{i,k} \times I_{j,k})$ попарно различны для различных троек $(i,j,k)$. 
    В самом деле, предположим, что $\tilde\varphi_k(I_{i,k} \times I_{j,k}) = \tilde\varphi_n(I_{p,n} \times I_{q,n})$. 
    Тогда $\varphi_k(I_{i,k}) = \varphi_n(I_{p,n})$. Поэтому $k=n$ и $i=p$. Аналогично, $j=q$.
        
    Возьмём кусочно-линейную функцию $h\colon[\,0,3\,]\to\R$ такую, что $h\bigl(\tilde\varphi_k(I_{i,k} \times I_{j,k})\bigr) = \frac{1}{3} f(l_{i,k}^0, l_{j,k}^0)$ для любых $k=\tofive$ и $i,j \in Z_k$. 
    Тогда $\norm{h} \leqslant \frac13\norm{f}$.
    Теперь лемма следует из того, что для любого $z \in I^2$
    \[
        \babs{S_\varphi h(z) - f(z)} \leqslant \frac{1}{6}\norm{f} + \frac{2}{3}\norm{f} < \frac{6}{7}\norm{f}.
    \]
    Здесь второе неравенство очевидно.
    Первое неравенство доказывается так.
    В силу выбора числа $N$ и равномерной непрерывности, $\abs{f(z) - f(z')} < \frac{1}{6}\norm{f}$ для любых $k=1,\ldots,5$, $i,j \in Z_k$ и $z,z' \in I_{i,k} \times I_{j,k}$.
    По утверждению \ref{propQuarters}, существует не менее трёх троек
    $(i_s,j_s,k_s) \in Z_{k_s} \times Z_{k_s} \times \tofiveset$, $s = 1,2,3$, таких, что $I_{i_s,k_s}\times I_{j_s,k_s} \ni z$.
    Для них $h\bigl(\tilde\varphi_{k_s}(z)\bigr) = \dfrac{1}{3}f(l_{i_s,k_s}^0, l_{j_s,k_s}^0)$. 
    Эти значения отличаются от $\dfrac{1}{3}f(z)$ не более, чем на $\dfrac{1}{18}\norm{f}$.
    Следовательно,
    \[
     \abs{\sum_{s=1}^3 h\bigl(\tilde\varphi_{k_s}(z)\bigr) - f\bigl(z\bigr)} \leqslant \dfrac{1}{6}\norm{f}.
    \]
    Другие два значения $h\bigl(\tilde\varphi_k(z)\bigr)$ не превосходят $\dfrac{1}{3}\norm{f}$ по абсолютной величине.
    \hfill$\square$
